\documentclass[10pt,twocolumn]{article} 
\usepackage{simpleConference}
\usepackage{times}
\usepackage{graphicx}
\usepackage{amssymb}
\usepackage{amsmath}
\usepackage{url,hyperref}
\usepackage{subcaption}
\usepackage{authblk}
\begin{document}

\title{Epidemic Conditions with Temporary Link Deactivation on a Network SIR Disease Model}

\author{Hannah Scanlon$^1$ and John Gemmer $^1$ \\
$^1$ Department of Mathematics and Statistics, Wake Forest University, Winston Salem NC 27109  \\
\vspace{15 pt}
\today \\
}


\maketitle
\thispagestyle{empty}

\begin{abstract}
The spread of an infectious disease depends on intrinsic properties of the disease as well as the connectivity and actions of the population. This study investigates the dynamics of an SIR type model which accounts for human tendency to avoid infection while also maintaining preexisting, interpersonal relationships. Specifically, we use a network model in which individuals probabilistically  deactivate connections to infected individuals and later reconnect to  the same individuals upon recovery. To analyze this network model, a mean field approximation consisting of a system of fourteen ordinary differential equations for the number of nodes and edges is developed. This system of equations is closed using a moment closure approximation for the number of triple links. By analyzing the differential equations, it is shown that, in addition to force of infection and recovery rate, the probability of deactivating edges and the average node degree of the underlying network determine if an epidemic occurs.
\end{abstract}



\section{Introduction}
The COVID-19 pandemic has had a profound impact on society. In response, the mathematics and broader scientific community has focused considerable research efforts to understand the spread of the virus and its impact not only on physical health \cite{Goldstein22035} but on mental health \cite{Witteveen202009609}, the economy~\cite{Witteveen202009609, bonaccorsi2020economic}, policy~\cite{flaxman2020estimating, qiu2020nationwide}, climate~\cite{Venter18984}, distribution networks~\cite{bonaccorsi2020economic}, equitable distribution of vaccines~\cite{bollyky2020equitable}, and racial disparities \cite{laurencin2020covid, chowkwanyun2020racial, Wrigley-Field21854} to name but a few. Despite a tremendous volume of research in this area, there
is still considerable effort devoted to developing and analyzing improved mathematical models that address aspects of the above issues. In particular, there is a clear need for epidemiological models that incorporate human behavior.

In this paper we propose and study a model for the spread of an infectious disease on an adaptive network in which individuals can temporarily deactivate connections with infected individuals and then reconnect upon recovery. Such a situation could arise, for example, in an office setting in which infected employees reduce their work hours or stay at home all together and thus lower their average number of contacts in a day. The problem we address in our model is the determination of a minimal deactivating rate needed to eliminate the spread of the disease as a function of the average node degree of the network, the force of infection, and the recovery rate of the disease. 

Naturally, in an adaptive network the spread of the disease can be eliminated by deleting or isolating all connections with infected individuals. However, for realistic human networks the implementation of such a process through stay at home orders or lockdown of businesses could be infeasible for a variety of reasons, e.g. the work force consists of essential workers, compliance may not be absolute, the economic impact would be too extreme \cite{bonaccorsi2020economic}, etc. Instead, by implementing an intermediary deactivation rate the network can still be productive since some connections are maintained while the spread of the disease is mitigated. 

\subsection{Background and drawbacks of classic models}
\begin{figure}
     \centering
     \begin{subfigure}[b]{0.2\textwidth}
         \centering
         \includegraphics[width=.8\textwidth]{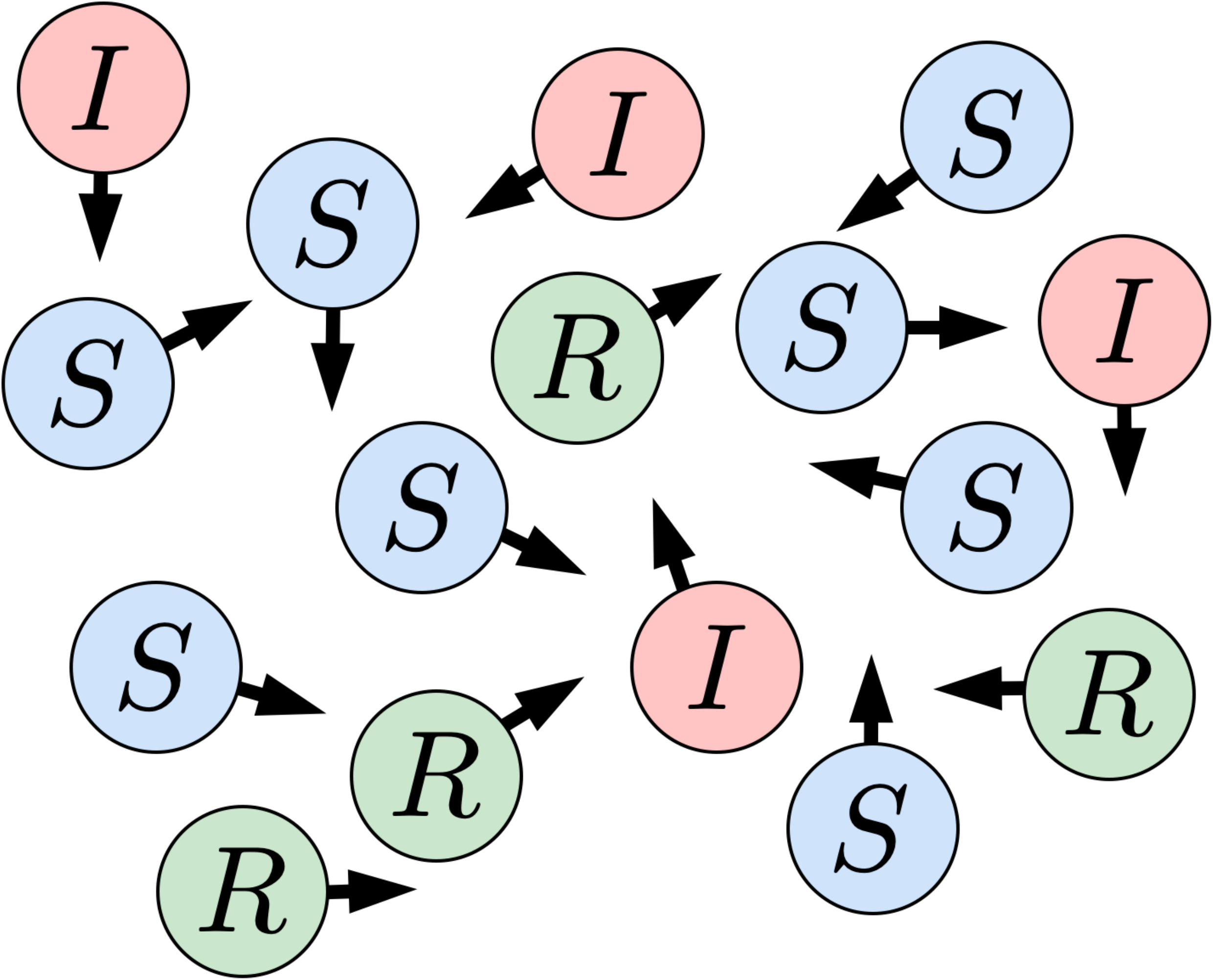}
         \caption{}
         \label{fig:micro}
     \end{subfigure}
     \hspace{.2in}
     \begin{subfigure}[b]{0.2\textwidth}
         \centering
         \includegraphics[width=.8\textwidth]{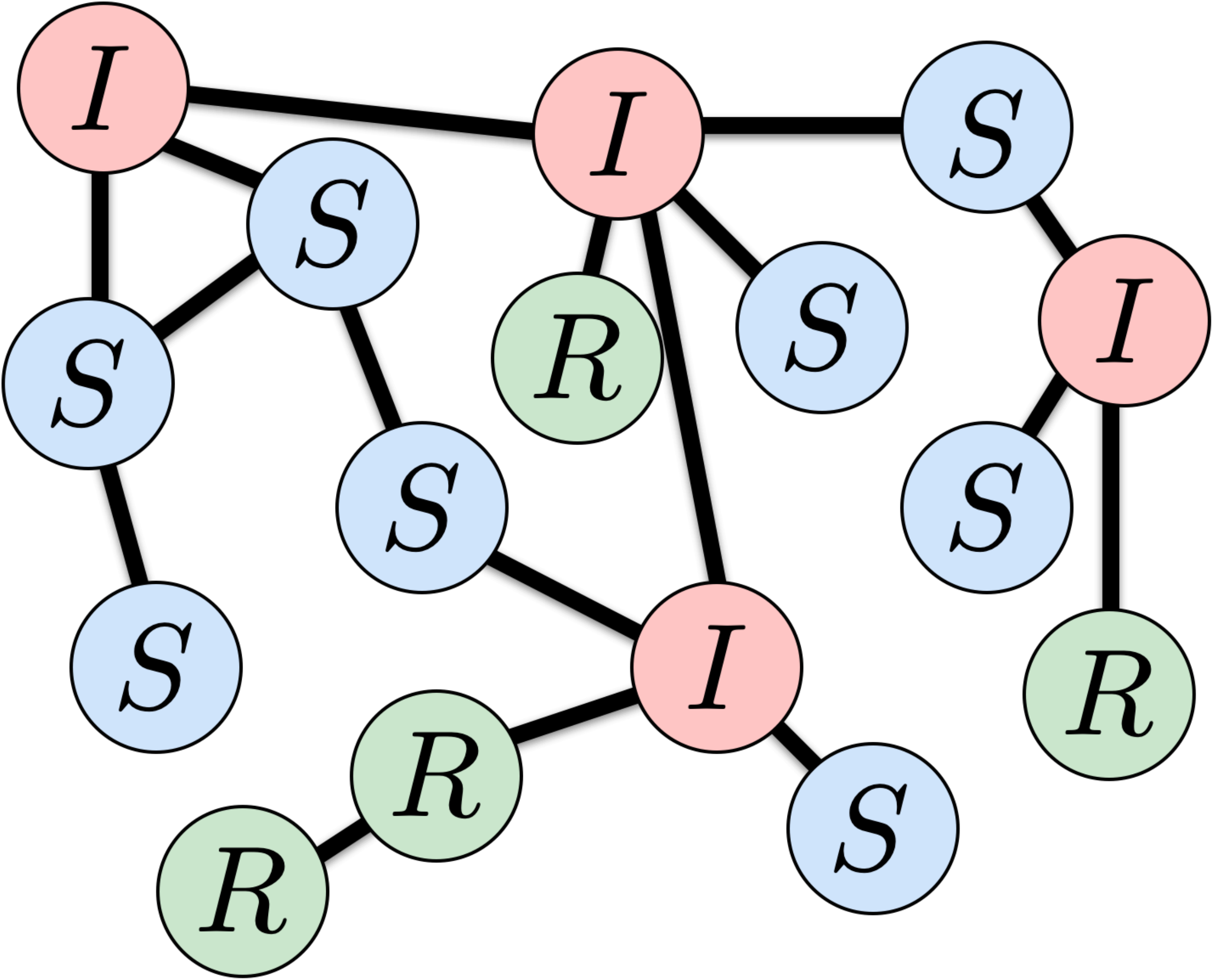}
         \caption{}
         \label{fig:meso}
     \end{subfigure}
     \begin{subfigure}[b]{0.2\textwidth}
         \centering
         \includegraphics[height=.7\textwidth]{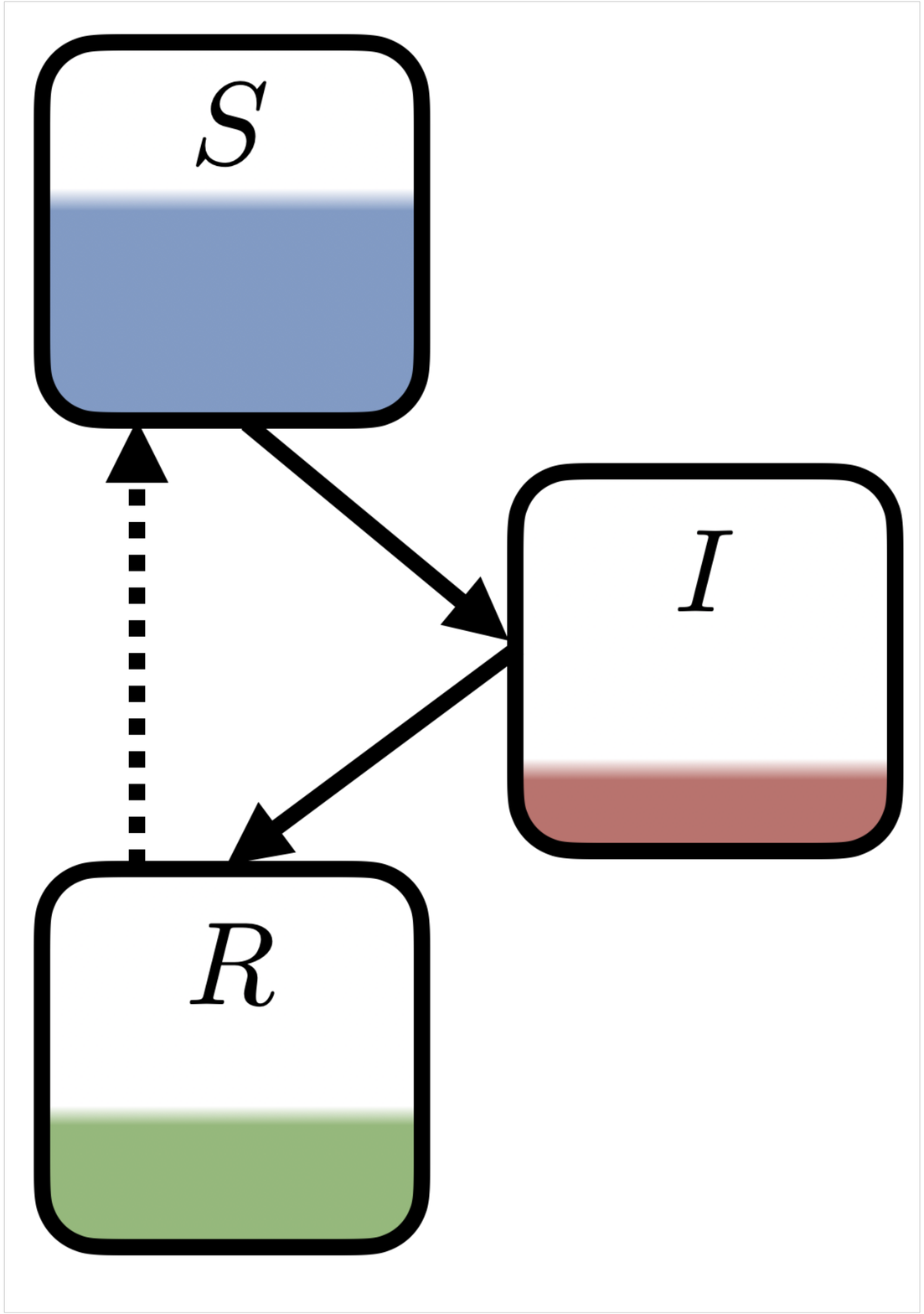}
         \caption{}
         \label{fig:macro}
     \end{subfigure}
      \begin{subfigure}[b]{0.2\textwidth}
         \centering
         \includegraphics[height=.7\textwidth]{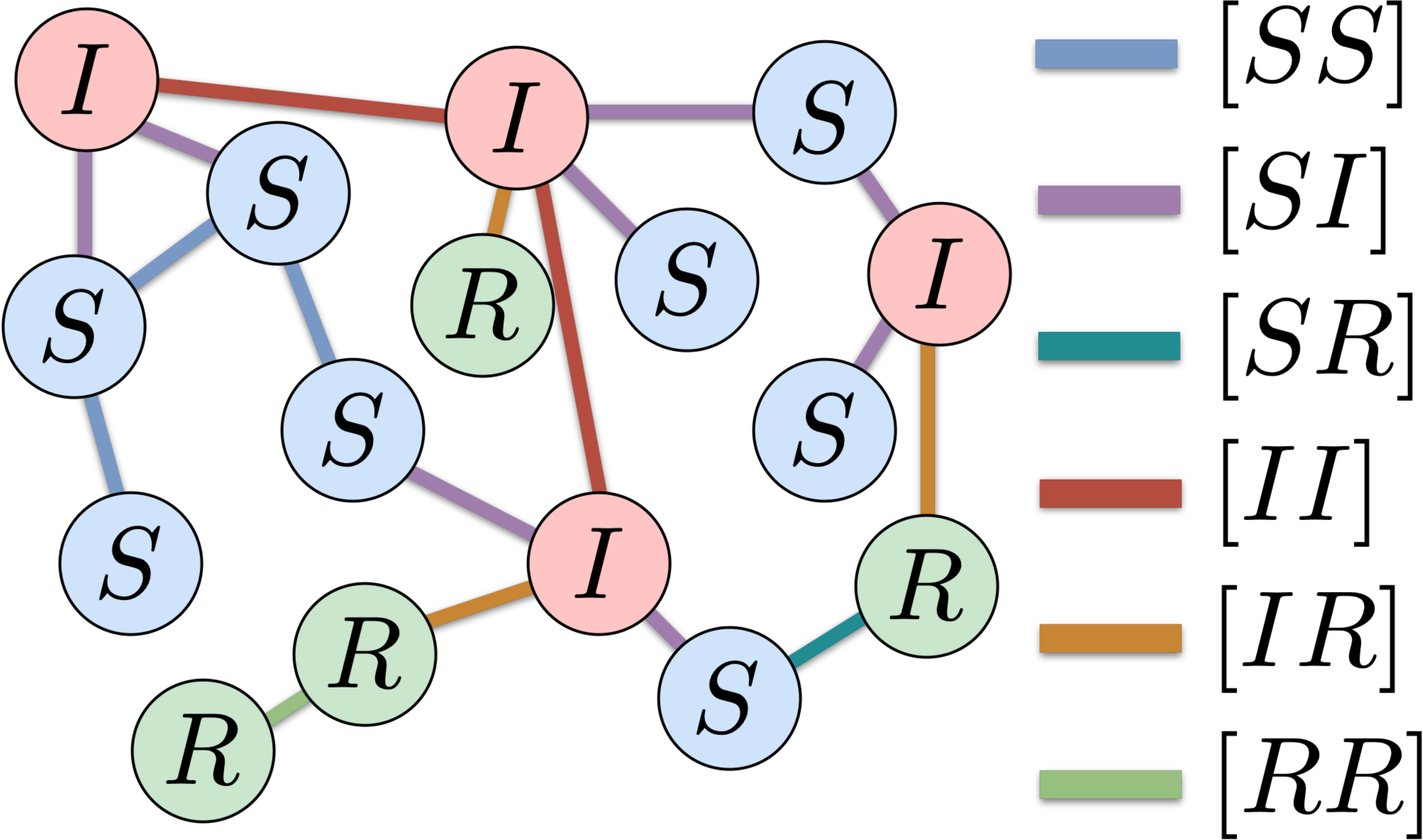}
         \caption{}
         \label{fig:adaptive}
     \end{subfigure}
        \caption{Illustration of mathematical frameworks for modeling the spread of infections diseases in which the status of each individual is either susceptible $(S)$, infected $(I)$, or recovered $(R)$. (a) Agent based model coupling spatial dynamics with the spread of the disease. Arrows indicate the direction of motion of each individual in the system. (b) Static network with the spread of the disease propagating along edges. (c) Compartment model with the disease spreading between the various population densities. (d) Adaptive network in which the population and edge densities are incorporated into compartment models.}
        \label{fig:MacroDescription}
\end{figure}

Before beginning a discussion of modeling the spread of infectious disease on adaptive networks, we first step back and discuss classic models for the spread of infectious diseases. There are a large number of mathematical models for the spread of infectious diseases whose efficacy and validity vary over a wide range of spatial and temporal scales. Typical mathematical models consist of agent based models at the microscale~\cite{rousseau1997dynamical, fuentes1999cellular, keeling2011modeling, tang2020review}, to network models at the mesoscale~\cite{Gross_2006, gross2009adaptive, nowzari2016analysis, Demirel_2017, kiss2017mathematics}, to finally mean-field compartment models at the macroscale~\cite{hethcote2000mathematics, allen2008mathematical, martcheva2015introduction,brauer2017mathematical, brauer2019mathematical}; see Figure~\ref{fig:MacroDescription}(a)--(c). Following the pioneering of Kermack and McKendrick~\cite{Kermack}, the unifying thread in all of these models is that members of the population are categorized depending on their infection status, e.g., susceptible ($S$), infected ($I$), and recovered ($R$), and the dynamic evolution of each individual's status is modeled either as a stochastic or purely deterministic process. In agent based models this consists of providing rules for the movement of individual agents as well as the transmission of infection between susceptible and infected agents. In network models infection between individuals occurs along undirected edges of an underlying static contact network. Finally, in compartment models the state variables consist of fractions of the total population with a given infection status and the disease evolves according to a differential equation. 

For reference, the standard compartment models are the $SIR$ model given by
\begin{equation}
    \begin{aligned}
        \dot{S}&=-\beta SI,\\
        \dot{I}&=\beta SI-\gamma I,\\
        \dot{R}&=\gamma I,
    \end{aligned}
\end{equation}
and the $SIS$ model given by
\begin{equation}
\begin{aligned}
    \dot{S}&=-\beta SI+\gamma I,\\
    \dot{I}&=\beta SI -\gamma I,
    \end{aligned}
\end{equation}
where $\beta$ is the per capita infection rate and $\gamma$ the recovery rate \cite{hethcote2000mathematics}. The $SIR$ model is often used to study the spread of diseases that confer lifelong immunity while the $SIS$ model is commonly used to study the spread of sexually transmitted diseases.

The network and compartment approaches can be linked by Kurtz's theorem which is  essentially a law of large numbers which states that the average dynamics of the Markov process at the network level limits to the dynamics of the deterministic differential equations at the compartment level as the size of the network $N\rightarrow \infty$ \cite{kurtz1971limit, nowzari2016analysis}. Such a limit we will refer to as a continuum limit.


The benefit of adopting a compartment modeling approach is that they are amenable to mathematical analysis since standard tools from dynamical systems such as bifurcation theory can be used to precisely quantify conditions under which the number of infected individuals grows in time. These conditions are often given in terms of the basic reproduction number $R_0>0$ which is the number of individuals  a single infected individual infects in a fully susceptible population \cite{hethcote2000mathematics,diekmann1990definition, van2002reproduction, martcheva2015introduction}. Specifically, when $R_0>1$, the infected population will grow causing an epidemic while if $R_0<1$ the disease will be eliminated. In the classic $SIR$ and $SIS$ compartment models $R_0$ can be explicitly calculated and is given by the ratio $R_0=\beta/\gamma$ \cite{hethcote2000mathematics}.

Classic compartment models such as $SIR$ and $SIS$ are useful models in predicting the spread of a disease on short timescales but they have a number of drawbacks that limit their efficacy on larger timescales. Namely, such models assume i) a constant population size, ii) a low number of states, iii) a well mixed population and iv) there is no feedback between human behavior and the spread of the disease \cite{nowzari2016analysis}. The first two drawbacks can be addressed by incorporating population growth into the classic models as well as introducing additional compartments, e.g. exposed $(E)$, treatment $(T)$, quarantine $(Q)$ and vaccinated $(V)$ \cite{van2002reproduction, kiss2010impact,hikal2014dynamic}. The third drawback is equivalent to the assumption that the underlying contact network is given by a complete graph. The fourth drawback can be heuristically addressed by introducing new compartments or by allowing parameters like the infection and recovery rates to depend on the state variables. However, problems with introducing a large number of new compartments include the system might become intractable to analysis and the introduction of a large number of parameters could obscure the physical mechanisms which govern the spread of the disease. Moreover, by allowing parameters like the infection rate in a compartment model to depend on the state of the disease, it is not clear that such a model could be obtained in a continuum limit from an underlying network model.  

\subsection{Background on adaptive networks}
With the prior discussion serving as a backdrop, we now discuss adaptive network models for the spread of infectious diseases that more naturally incorporate human behavior; see for instance \cite{Gross_2006, Shaw_2008, marceau2010adaptive, schwartz2010rewiring, network_vaccine, TemporaryLinkDeactivation1, TemporaryLinkDeactivation2, schwartz2010rewiring, network_vaccine}. The key idea in such models is that individuals can change the topology of the network depending on the infection status of their contacts. For example, susceptible individuals could replace contacts with infected individuals with connections to susceptible individuals as in \cite{Gross_2006,  marceau2010adaptive}, or delete contacts with infected individuals as in \cite{TemporaryLinkDeactivation1, TemporaryLinkDeactivation2}. In this framework, in addition to the infection status of individual nodes, the edges themselves are also given a status depending on the infection status of the nodes connected by the edge. For example, for an $SIS$ model on a network the three states of the edges are given by $[SS]$, $[SI]$ and $[II]$ denoting the status of an edge connecting two susceptible nodes, an infected and susceptible node, and two infected nodes respectively. Figure \ref{fig:MacroDescription}(d) illustrates the nine resulting state variables for an $SIR$ model on a network placed within this framework. The adaptive network model then typically assumes that edges with an infected component, i.e. an $[SI]$ edge, will change its status with some probability to reduce the spread of the disease amongst nodes, e.g. an $[SI]$ edge rewires to different nodes to create an $[SS]$ edge with some probability.  

The average dynamics on an adaptive network can also be approximated by the dynamics of a compartment model in an appropriate continuum limit. For example, on a static network, the governing equations for an $SIS$ model incorporating edge dynamics is given by:
\begin{equation}
    \begin{aligned}
        \dot{S}&=-\beta [SI]+\gamma I,\\
        \dot{I}&=\beta [SI]-\gamma I,\\
        \dot{[SS]}&=\gamma[SI]-\beta [SSI],\\
       \dot{[SI]}&=\beta\left([SSI]-[SI]-[ISI]\right)-\gamma\left([SI]-2[II]\right),\\
        \dot{[II]}&=\beta\left([SI]+[ISI]\right)-2\gamma [II],
    \end{aligned} \label{Eqn:SIS}
\end{equation}
where $[ABC]$ denotes the density of triple links with a given sequence of states $A,B,C\in \{S, I\}$ \cite{Gross_2006, kiss2017mathematics}. The first two equations model the infection and recovery of nodes while the remaining equations correspond to the conversion of various edge types as nodes are infected or recover. The state variables in the above equations are implicitly understood to correspond to the expected values of the node and edge densities, however the notation $\mathbb{E}$ for expectation is suppressed. If we assume further that $\mathbb{E}\left([SI]\right)=\mathbb{E}(S)\mathbb{E}(I)$, i.e. assume a well mixed population, we obtain the standard $SIS$ model. However, the benefit of retaining the dynamics of the edges is that human behavior can now be incorporated  directly into Equation \eqref{Eqn:SIS} by modifying  the dynamics of $[SS]$, $[SI]$, and $[II]$ while retaining the same dynamics on $S$ and $I$. 

The drawback of the continuum limit presented in Equation \eqref{Eqn:SIS} is that it does not form a system of closed equations. Specifically, the dynamics of the triple links must be specified resulting in the need for equations governing the quartic links and so on. In order to close the system at the level of the dynamics for the edges, the number of triple links must be approximated by using a process called a moment closure. The simplest moment closure can be derived by assuming a homogeneous degree distribution and applying a counting argument. This moment closure is given by:
\begin{equation}
    [ABC]\approx \frac{\langle k \rangle-1}{\langle k \rangle }\frac{[AB][BC]}{B}, \label{Eqn:MomentClosure}
\end{equation}
where $\langle k \rangle$ is the average degree of a node \cite{kiss2017mathematics}; see the Appendix for a derivation. More sophisticated moment closures that account for inhomogeneities in the degree distribution arising from the friendship paradox, existence of triangles, a high clustering coefficient, etc. can be derived based on the topology of the network; see for instance \cite{TemporaryLinkDeactivation1, TemporaryLinkDeactivation2, marceau2010adaptive}.

\subsection{A roadmap}
We conclude the Introduction with a roadmap for the paper. In Section 2 we present mathematical models for the spread of an infectious disease on adaptive small world networks at both the network and compartment level. Our models are built on the work of Shaw et. al. in \cite{TemporaryLinkDeactivation1, TemporaryLinkDeactivation2} in which an SIS model was implemented on an adaptive network with temporary link deactivation. Network assumptions such as those used by Gross et al. in \cite{Gross_2006}, prioritize maintaining the original connectedness, or average node degree, of a network and reflects some aspects of human interactions by disconnecting potentially infectious connections and creating new, safer links. On an interpersonal scale, however, we know this to be inconsistent with human behavior. To address this concern, our adaptive network model preserves known relationships throughout the course of the disease while allowing individuals to protect themselves from infection by temporarily deactivating potentially infectious interactions. 

In Section 3 we present the primary results of our work. We first numerically study the convergence between the network and compartment models in the continuum limit. While the compartment model slightly overestimates the dynamics of the disease, the edge dynamics agree remarkably well and moreover the parameter conditions in which an epidemic occur are in agreement. To further probe the conditions under which an epidemic occurs, we investigate parameter regimes in our compartment model in which state changes including not only the sign of $\dot{I}(0)$ but also $\ddot{I}(0)$ and $\ddot{S}(0)$ occur. Through these calculations we replicate the standard value of $R_0$ as well as identify the following critical edge deactivation rates: 
\begin{equation}
\begin{aligned}
    p_1^*&=\beta \left(\frac{\langle k\rangle}{2}-\frac{3}{2}\right)-\gamma,\\
    p_2^*&=p_1^*-\gamma +\frac{\gamma^2}{\beta \langle k \rangle}.
    \end{aligned}\label{eqn:critica_rates}
\end{equation}
Specifically, if the deactivation rate $p$ satisfies $p>p_1$ then $\lim_{I(0)\rightarrow 0}I(0)^{-1}\ddot{S}(0)>0$ and if $p>p_2^*$ then $\lim_{I(0)\rightarrow 0}I(0)^{-1}\ddot{I}(0)<0$. We provide numerical evidence that if $R_0>1$ but the deactivating rate is above these thresholds then the the disease will still initially spread but the total number of infected individuals is drastically reduced. This provides additional criteria beyond $R_0<1$ for controlling the spread of an infectious disease.  

We conclude in Section 4 with a discussion of our key results, the implications of our results that elucidate the connection between human behavior and the spread of a disease, and avenues for further work.


\section{Models}

This work uses two models to investigate temporary link deactivation on an $SIR$ disease model. The first is a network model which applies system changes including infection, recovery and edge deactivation as probabilities while tracking the states of all individual nodes and edges. We determine appropriate values for the number of Monte-Carlo simulations $M$, the temporal spacing $\Delta t$, and network size $N$ to ensure convergence of our simulations for the mean field dynamics. Using the determined parameter values, we consider the network model to be a proxy for reality since all dynamics are tracked on an individual scale. The second model is an ODE model with compartments for each node and edge type. This model approximates the network behavior on a macro scale while applying system changes as rates applied to the compartments. By developing an ODE model that reflects the network model behavior, we can more efficiently simulate and more robustly analyze the system behavior.

\subsection{Network Model}

For the network model, we study disease spread on a population represented by a graph, $G = \{\mathbf{V},\mathbf{E}\}$, where $\mathbf{V}$ denotes the set of $N$ vertices (i.e. nodes) and $\mathbf{E}$ denotes the set of edges. The graph used is a Watts-Strogatz model which creates a realistic model of human connections referred to as a small-world network \cite{WattsNetwork}. Specifically, this graph has a large number of nodes, short average path lengths and tightly knit groups of nodes as measured by a high clustering coefficient. This graph is created by first generating a  ring lattice of average node degree  $\langle k \rangle$. Applying the handshaking theorem with a constant degree yields the approximation $\bar{N} \approx \frac{\langle k \rangle N}{2}$ where $\bar{N}$ is the total number of edges in our system. A portion, $\alpha$, of the edges are then randomly rewired. This preserves the average node degree and total number of edges but creates the desired characteristics of a small world network including more tightly clustered nodes. From the graph, we generate an adjacency matrix which is a symmetric, $N\times N$ matrix, $A$, defined by $A_{i,j} = 1$ if node $i$ is connected to node $j$ and is $0$ otherwise. 

We model an SIR type disease progression on this network in which individuals move from susceptible ($S$) to infected ($I$) to recovered ($R$) corresponding to three possible node states $\{S, I, R\}$. Letting $i$ index nodes and $k$ index time, we define $V_i^k\in \{S, I, R\}$ as the state of node $i$ at time $k\Delta t$ where $\Delta t > 0$ is the temporal spacing. Based off of the status of $V_i^k$ we define another set of vectors $S_i^k, I_i^k, R_i^k$ with the $i^{\text{th}}$ entry equal to $1$ if $V_i^k = S,I,R$ respectively and $0$ otherwise. We define $\beta \Delta t$ as the probability of infection applied based on edges between a susceptible and infected node and $\gamma \Delta t$ as the probability of recovery applied to infected nodes. Finally, we apply a temporary deactivation assumption to the edges of the graph by storing deactivated edges in another symmetric, $N\times N$ adjacency matrix, $D$, which is initialized with all zeros.

The probability a node satisfying $V_{i}^{k} = S$ becomes infected at the next time step is equal to $\beta \Delta t$ times the number of active connections between that susceptible node and other infected nodes. The number of such connections is found by taking the difference between the adjacency matrix, $A$, and the current deactivated matrix, $D$, isolating the $i^{\text{th}}$ node's connections by right multiplying by the standard basis vector $e_i$ and summing the number of infected connections by left multiplying by the transpose of the vector $I^k$. The probability that a susceptible node remains susceptible is then $1$ minus the above calculated probability. The probability a node satisfying  $V_{i}^{k} = I$ recovers at the next time step is equal to $\gamma \Delta t$ while the probability that the same node remains infected is $1- \gamma \Delta t$. Finally, the probability a node satisfying $V_{i}^{k} = R$ remains recovered is equal to 1 since we assume the recovered class is immune and cannot return to susceptible or infected states. This gives the infection probabilities 

\begin{equation}
\begin{aligned}
    \mathbb{P}(V_{i}^{k+1} &= S | V_i^k = S) = 1-\beta\Delta t\cdot (I^k)^T(A-D)e_i, \\
    \mathbb{P}(V_{i}^{k+1} &= I | V_i^k = S) = \beta\Delta t\cdot (I^k)^T(A-D)e_i, \\
    \mathbb{P}(V_{i}^{k+1} &= I | V_i^k = I) = 1-\gamma\Delta t,\\
    \mathbb{P}(V_{i}^{k+1} &= R | V_i^k = I) = \gamma\Delta t,\\
    \mathbb{P}(V_{i}^{k+1} &= R | V_i^k = R) = 1.
\end{aligned}
\end{equation}

For the link deactivation assumption, since infection can be passed by connections between infected and susceptible nodes, we apply a temporary deactivation probability, $p \Delta t$, to any such edges. That is, if at time $k \Delta t$, $A_{i,j} = 1$, $D_{i,j}^k = 0$, $S_i^k = 1$, and $I_j^k = 1$ then with probability $p \Delta t$, $D_{i,j}^{k+1}$ and $D_{j,1}^{k+1}$ become $1$. When deactivated edges $D_{i,j}^k = 1$ are no longer potentially infectious, i.e. $i,j$ indices correspond to susceptible to recovered and recovered to recovered edges, they are reactivated with probability $r \Delta t$. That is, if at time $k \Delta t$, $D_{i,j}^k = 1$, $R_i^k = 1$ and $S_j^k = 1$ then with probability $r \Delta t$, $D_{i,j}^{k+1}$ and $D_{j,i}^{k+1}$ become $0$. Similarly, if at time $k \Delta t$, $D_{i,j}^k = 1$, $R_i^k = 1$ and $R_j^k = 1$ then with probability $\frac{r}{2} \Delta t$, $D_{i,j}^{k+1}$ and $D_{j,i}^{k+1}$ become $0$. Note this probability is halved to account for the symmetry of $R$ to $R$ edges. These assumptions preserve the original graph structure as created by the Watts-Strogatz model by never updating the adjacency matrix $A$. This gives the edge transition probabilities

\begin{equation}
\begin{aligned}
    \mathbb{P}(D_{i,j}^{k+1} &= 1 | S_i^k = I_j^k = A_{i,j} = 1, D_{i,j}^k =0) = p \Delta t, \\
    \mathbb{P}(D_{i,j}^{k+1} &= 0 | R_i^k = S_j^k = D_{i,j}^k =1) = r\Delta t, \\
        \mathbb{P}(D_{i,j}^{k+1} &= 0 | R_i^k = R_j^k = D_{i,j}^k =1) = \frac{r}{2}\Delta t,
\end{aligned}
\end{equation}
where transitions made to $D_{i,j}^{k+1}$ are made symmetrically to  $D_{j,i}^{k+1}$ for all $i$ and $j$.

\begin{figure}[htp]
     \centering
     \begin{subfigure}[b]{.2\textwidth}
         \centering
         \includegraphics[width = \linewidth]{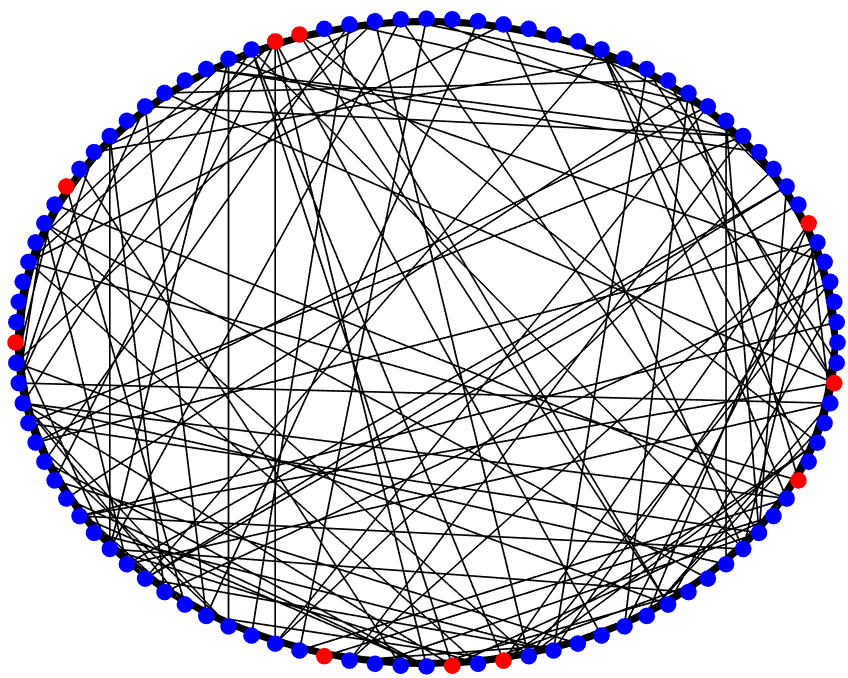}
         \caption{Initial State}
         \label{fig:InitialNetwork}
     \end{subfigure}
     \hfill
     \begin{subfigure}[b]{.2\textwidth}
         \centering
        \includegraphics[width =
    \linewidth]{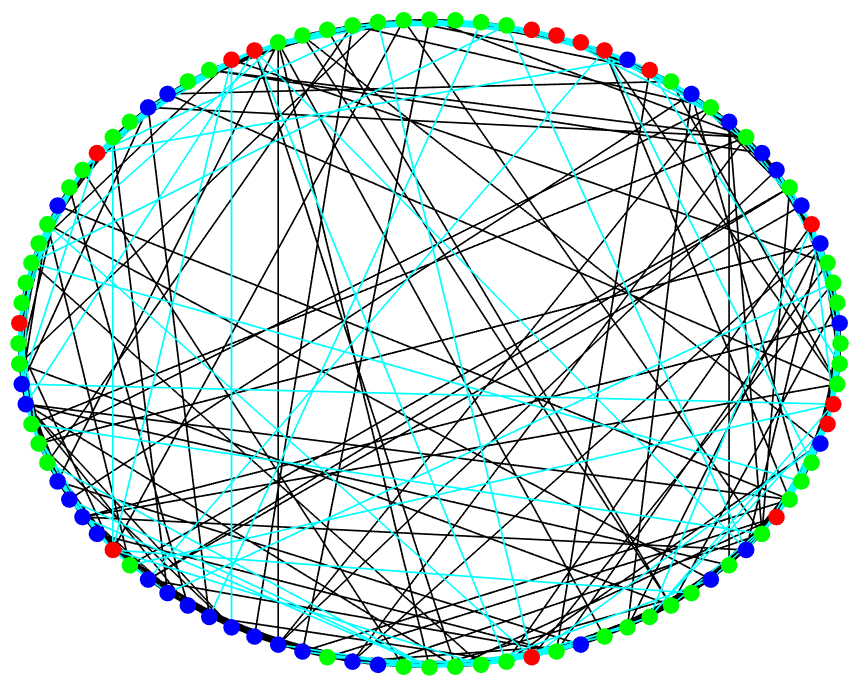}
         \caption{Intermediate State}
         \label{fig:MidRunNetwork}
     \end{subfigure}
     \hfill
     \begin{subfigure}[b]{.2\textwidth}
         \centering
        \includegraphics[width =
    \linewidth]{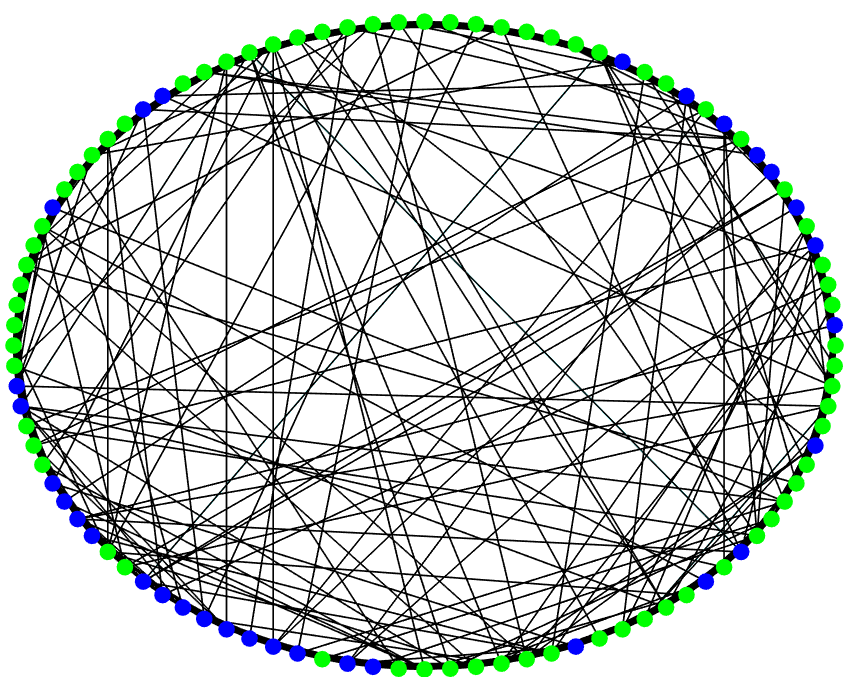}
         \caption{Final State}
         \label{fig:EndNetwork}
     \end{subfigure}
        \caption{Snapshots of the dynamics of the $SIR$ network model with link deactivation at (a) time $t=0$, (b) an intermediate time, and (c) the final network state. Black lines correspond to active edges while cyan are temporarily deactivated edges. Susceptible nodes are blue, infected red, and recovered green. The initial conditions  consisted of $10$ randomly selected infected nodes and the parameters were given by $\beta=0.1$, $\gamma=0.2$, $p=0.8$, and $r=0.9$.}
    \label{fig:NetworkProgression}
\end{figure}

Figure \ref{fig:NetworkProgression} depicts an example progression of the Watts Strogatz network states through a disease simulation before, during, and after the infectious event. In all simulations, the network was initialized with 10\% of nodes randomly selected to be infected.  

\subsection{Convergence of Network Model}

To ensure the consistency of conclusions drawn from the statistics of Monte-Carlo (MC) simulations of our network model we need to test for convergence in the number of simulations $M$, time step $\Delta t$ and network size $N$. For all convergence analysis, we will use the $L_2$ norm as our diagnostic for convergence. Specifically, to compute an error for $M$ simulations we will generate two sets of data of $M$ simulations each assuming one of these sets is a proxy for the converged statistics. On each set of $M$ simulations we compute the average number of infected nodes at time step $k$ and denote these computed values by $J^k$ and $\overline{J}^k$ respectively. The relative error is then computed using the $L_2$ norm and is given by
\begin{equation}
    E(M,N,\Delta t, p)=\frac{\left(\sum_{k}\left(J^k-\bar{J}^k\right)^2\right)^{\frac{1}{2}}}{\left(\sum_{k}\left(J^k\right)^2\right)^{\frac{1}{2}}},
\end{equation}
where we have also included $p$ as a variable to emphasize that the deactivating rate could influence the convergence. Note, this definition of the error is equivalent to estimating the variance of the Monte-Carlo estimator \cite{rubinstein2016simulation}. 

\begin{figure}[htp]
     \centering
     \begin{subfigure}[b]{.9\linewidth}
         \centering
         \includegraphics[width = \linewidth]{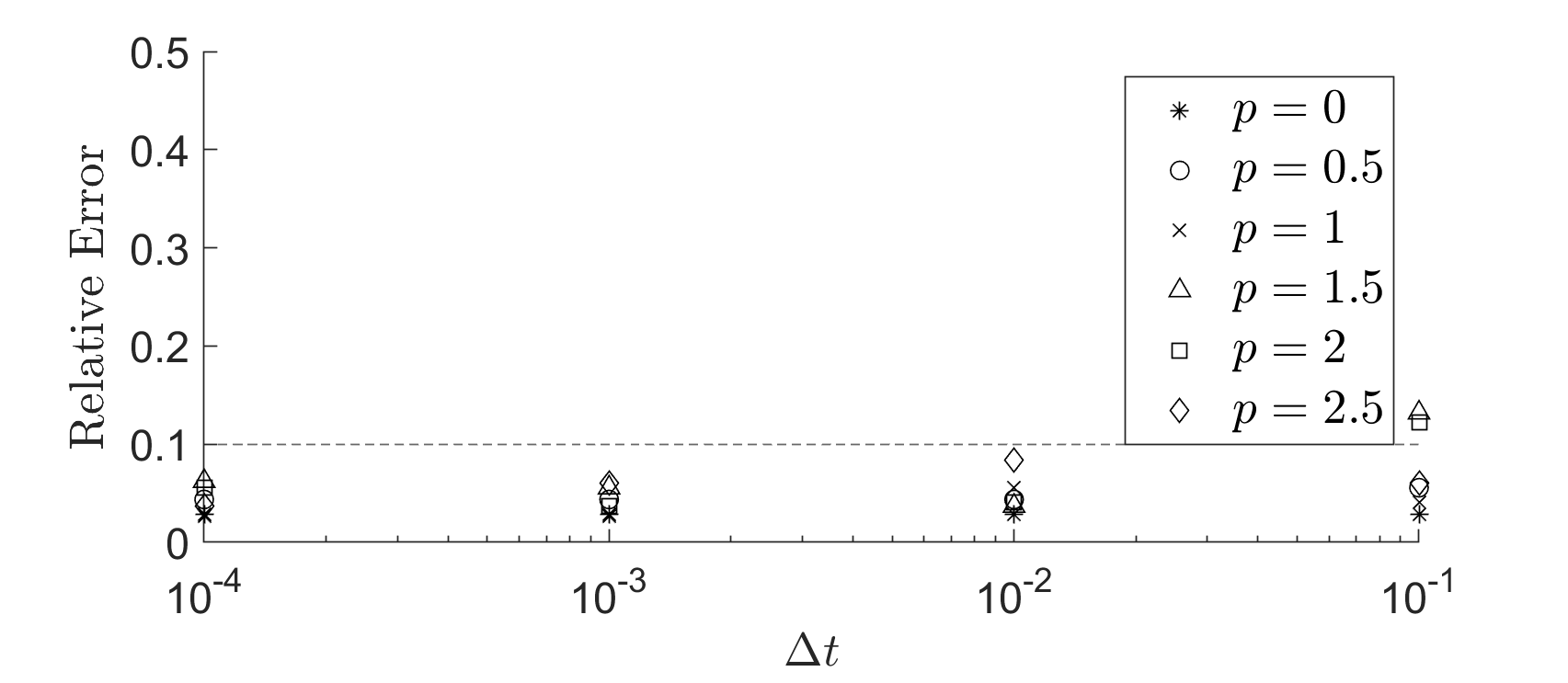}
         \caption{$100$ Node Convergence}
         \label{fig:dtMCerror}
     \end{subfigure}
     \hfill
     \begin{subfigure}[b]{.9\linewidth}
         \centering
        \includegraphics[width =
    \linewidth]{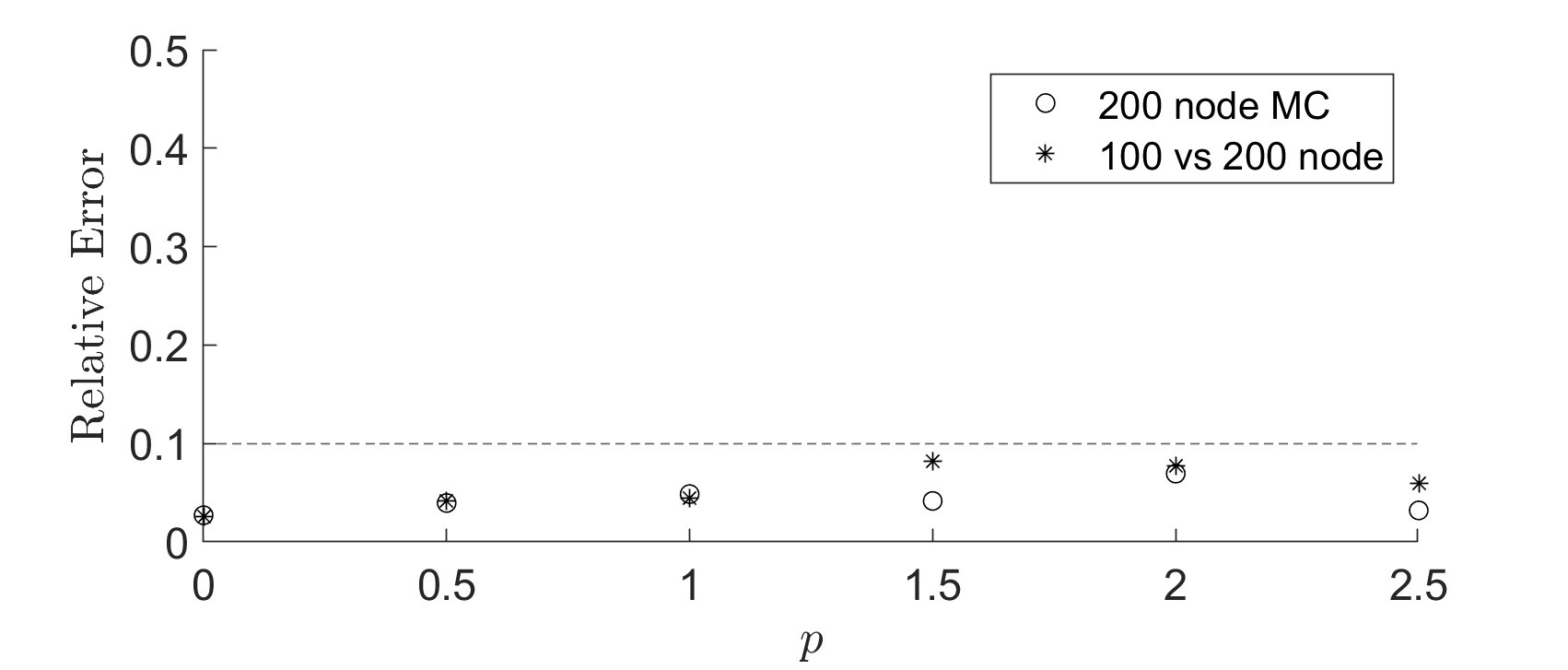}
         \caption{$200$ Node Convergence}
         \label{fig:sizeerror}
     \end{subfigure}
        \caption{Numerically computed error of the average number of infected individuals for MC simulations of the network model. (a) $100$ node network with $20$ MC simulations for various values of $\Delta t$ and $p$. (b) 200 node network with 2 sets of 20 MC simulations (circles) and 200 vs 100 node network (asterisk) for various $p$ and $\Delta t = 0.01$ }
    \label{fig:ErrorPlots}
\end{figure}

First, we investigate the convergence on a $100$ node network over a range of deactivation rates and time step sizes between two sets of $20$ MC simulations. Specifically, we compute $E(20,100,\Delta t, p)$ for $\Delta t$ values from $0.1$ to $0.0001$ and $p$ values ranging from $0$ to $2.5$. Figure \ref{fig:ErrorPlots}(a) shows the value of of our error is below our cut off value of $0.1$ for $\Delta t = 0.01$ and all tested $p$ values. Additionally, the error did not decrease significantly for smaller $\Delta t$ values. This indicates that averaging $20$ simulations with $100$ nodes and $\Delta t = 0.01$ gives sufficient convergence for any $p \in [0, 2.5]$.

Next, in order to compare the network model to a compartment model, we consider the continuum limit of the system as the number of nodes $N \to \infty$. We compute $E(20,200,0.01, p)$ for $p$ values ranging from $0$ to $2.5$. The results shown as circles in Figure \ref{fig:ErrorPlots}(b) shows $20$ MC simulations is also sufficient for the $200$ node network to meet our $0.1$ cut off for the relative error. 

For our final convergence analysis, we need a different definition of error. This error will compare the results between a set of simulations with 100 nodes and a set with 200 nodes. We compute the average number of infected nodes from each set of simulations at time step $k$ and denote these computed values from the 100 node network by $J^k$ and from the 200 node network by $L^k$. This relative error is again computed using the $L_2$ norm and is given by
\begin{equation}
    F(M,N_1,N_2,\Delta t, p)=\frac{\left(\sum_{k}\left(J^k-L^k\right)^2\right)^{\frac{1}{2}}}{\left(\sum_{k}\left(J^k\right)^2\right)^{\frac{1}{2}}}.
\end{equation}

We compute $F(20,100,200,0.01,p)$ for $p$ values ranging from $0$ to $2.5$. The results shown as asterisks in Figure \ref{fig:ErrorPlots}(b) demonstrate convergence of the results from the $200$ to the $100$ node network. Altogether, this analysis gives us confidence that our $100$ node network simulated with time step $\Delta t = 0.01$ for $p \in [0, 2.5]$ has converged sufficiently to compare results to a compartment ODE model. Table \ref{tab:NetworkParameters} summarizes the parameters previously defined and provides the values used in this study.

\begin{table}[t]
    \centering
    \renewcommand{\arraystretch}{1.2}
    \begin{tabular}{ |c|c|c| } 
 \hline
 Parameter & Definition & Value  \\ 
 \hline
  $\alpha$ & Watts-Strogatz rewiring & 0.2 \\ 
 \hline
  $\langle k \rangle$ & average node degree & 12 \\ 
 \hline 
 $\beta$ & infection probability & [0, 1]  \\ 
 \hline
$\gamma$ & recovery probability & 0.2  \\ 
 \hline
 $p$ & deactivating probability & [0, 2.5] \\ 
 \hline
 $r$ & reconnecting probability & 0.9 \\ 
 \hline
  $\Delta t$ & temporal spacing & .01 \\ 
 \hline
 $N$ & number of nodes & 100 \\ 
 \hline
  $\bar{N}$ & number of edges & 600 \\ 
 \hline
  $M$ & number of MC simulations & 20 \\ 
 \hline
\end{tabular}
    \caption{Parameter values used in MC simulations of the network model.}
    \label{tab:NetworkParameters}
\end{table}

\subsection{Proposed ODE Model}

In order to facilitate system analysis, we approximate the network with a system of ordinary differential equations using a mean-field approach. The following system of differential equations describes the change in the number of each type of node:

\begin{equation}
\begin{aligned}
    \dot{S} =& -\beta [SI],\\
    \dot{I} =& \beta[SI]-\gamma I, \\
    \dot{R} =& \gamma I, 
    \label{eq:ODENodes}
\end{aligned}
\end{equation}
and the number of each type of edge:

\begin{equation}
\begin{aligned}
    \dot{[SS]} =& -\beta[SSI],\\
    \dot{[SI]} =& \beta[SSI] - \beta([SI]+[ISI])-\gamma[SI] - p[SI],\\
    \dot{[SR]} =& -\beta[ISR] +\gamma[SI] + r[\widehat{SR}],\\
    \dot{[II]} =& \beta([SI]+[ISI])-2\gamma[II],\\
    \dot{[IR]} =& 2\gamma[II]-\gamma[IR]+\beta[ISR],\\
    \dot{[RR]}=& \gamma[IR]+r[\widehat{RR}],\\
     \dot{[\widehat{SI}]} =& p[SI]-\gamma[\widehat{SI}]-\beta[I\widehat{SI}],\\
    \dot{[\widehat{SR}]} =& \gamma[\widehat{SI}]-r[\widehat{SR}]-\beta[I\widehat{SR}],\\
    \dot{[\widehat{II}]} =& \beta[I\widehat{SI}]-2\gamma[\widehat{II}],\\
    \dot{[\widehat{IR}]} =& 2\gamma[\widehat{II}]+\beta[I\widehat{SR}]-\gamma[\widehat{IR}],\\
    \dot{[\widehat{RR}]}=& \gamma[\widehat{IR}]-r[\widehat{RR}], \label{eq:ODEEdges}
\end{aligned}
\end{equation}
where $X$ denotes the expected number of nodes of each type and $[XY]$ and $[\widehat{XY}]$ denote the number of active and deactivated edges respectively between nodes in state $X$ and $Y$ with $X \in \{S, I, R\}$. The notation $[X\widehat{YZ}]$ represents triple connection between an $[XY]$ edge and an $[\widehat{YZ}]$ edge centered at a $Y$ type node. For ease of presentation we wrote the above equations with the triple link states but these equations were closed using Equation \eqref{Eqn:MomentClosure}.

\begin{figure}
    \centering
    \includegraphics[width =.9
    \linewidth]{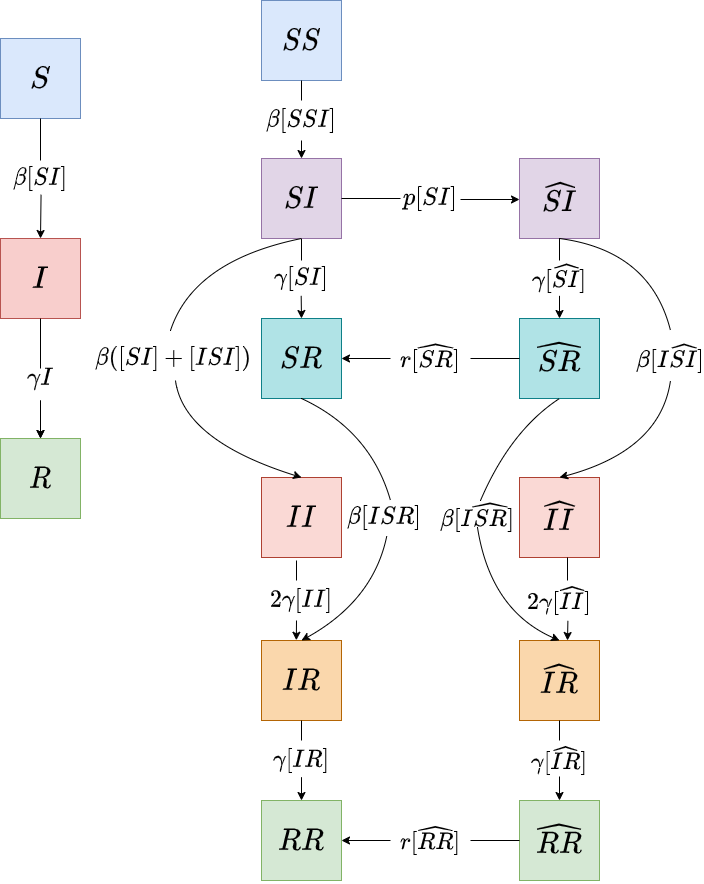}
    \caption{Flow chart depicting edge and node dynamics in the ODE model.}
    \label{fig:odeflow}
\end{figure}

The transitions in the ODE model reflect the behavior of the previously defined network model and utilize the same parameters. The state transitions are depicted in Figure \ref{fig:odeflow} with the node states, active edge states, and deactivated edge states shown in the first, second and third columns respectively. The node state transitions include susceptible nodes being infected at rate $\beta$ proportional to $[SI]$ edges and $I$ nodes recovering at rate $\gamma$. Edge state transitions involve the parameters $p$, $r$, $\beta$, and $\gamma$. The deactivation parameter $p$ is applied only to $[SI]$ which become $[\widehat{SI}]$. The reconnecting parameter $r$ is applied to both $[\widehat{SR}]$ and $[\widehat{RR}]$ which return to their equivalent active edge compartments. Edge states involving the infection of an $S$ node transition at rate $\beta$. For active edges, these transitions include $[SS]$ becoming $[SI]$ through $[SSI]$ triples, $[SI]$ becoming $[II]$ through both $[SI]$ and $[ISI]$ triples, and $[SR]$ becoming $[IR]$ through $[ISR]$ triples. Deactivated edges involve some of the same transitions including $[\widehat{SI}]$ becoming $[\widehat{II}]$ through $[I\widehat{SI}]$ triples and $[\widehat{SR}]$ becoming $[\widehat{IR}]$ through $[I\widehat{SR}]$ triples. Notably, the deactivated edges do not include a parallel $[SS]$ transition because $[\widehat{SS}]$ does not exist nor do they include $[\widehat{SI}]$  becoming $[\widehat{II}]$ directly through $[\widehat{SI}]$ since these deactivated edges cannot pass infection. Finally, the recovery rate $\gamma$ facilitates the transition of edges involved in the recovery of an $I$ node. The transitions include $[SI]$ becoming $[SR]$, $[II]$ becoming $[IR]$ at twice the recovery rate for each $I$ node involved, and $[IR]$ becoming $[RR]$. The same transitions occur in the equivalent deactivated edge compartments.

For all later simulations, we use ode45 in Matlab\cite{Matlab}. These equations become unstable as $S \to 0$ since the moment closure approximation divides by $S$. To account for this instability in our simulations, we set all moment closure approximations equal to $0$ when $S<0.001$. The system is initialized with $S =90$, $I = 10$ and $R = 0$ for the nodes. For the edges, we averaged the initial number of each edge type from 100 network simulations giving initial values $[SS] = 485$, $[SI] = 110$, $[II] = 5$ and $0$ for everything else. These initial conditions preserved the total number of nodes, edges and average node degree used in the network simulations. 

\section{Results}

\subsection{Infected Population Convergence}

\begin{figure}[htp]
     \centering
     \begin{subfigure}[b]{\linewidth}
         \centering
         \includegraphics[width = .82 \linewidth]{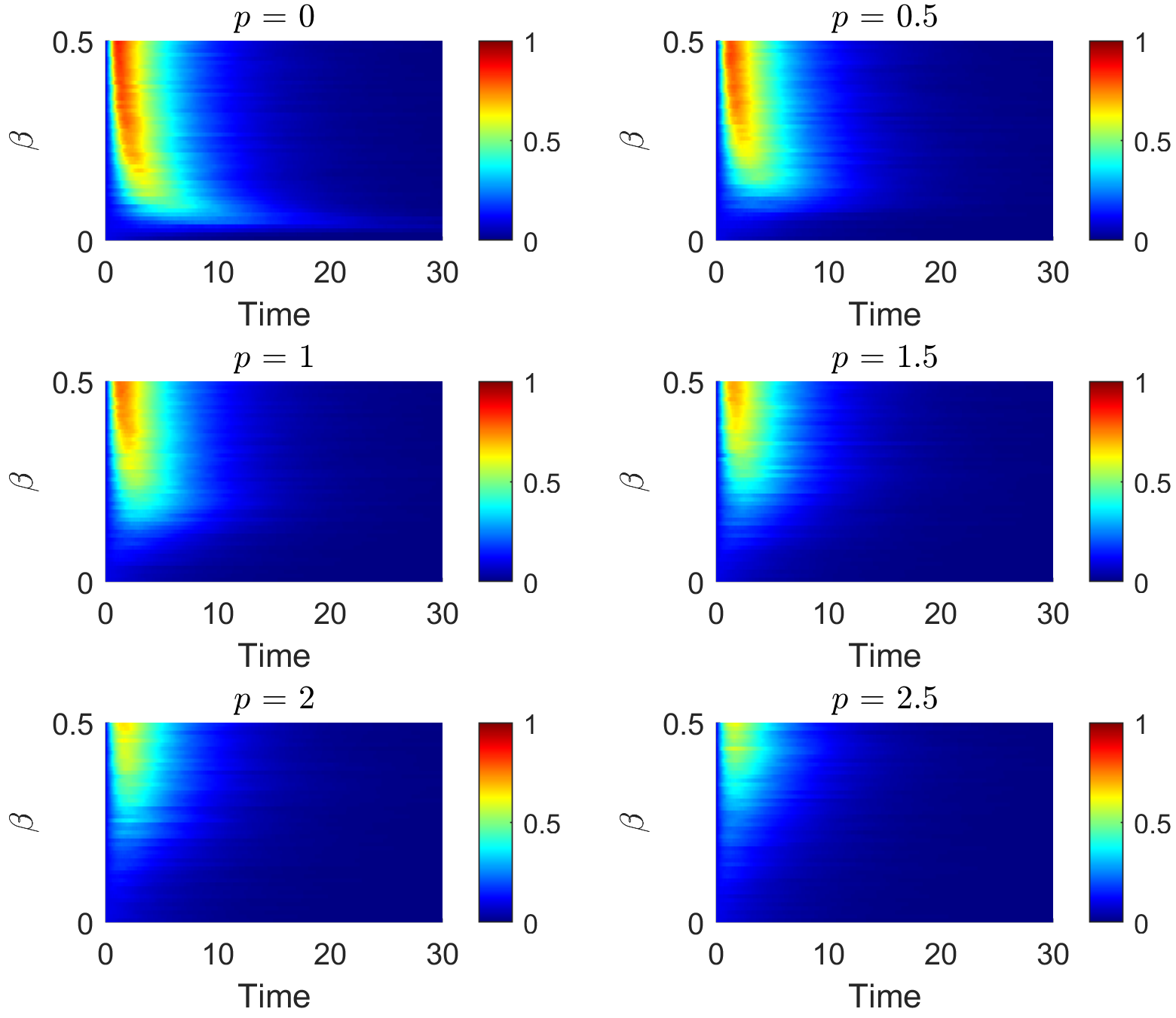}
         \caption{Network Model}
    \label{fig:NetworkContour}
     \end{subfigure}
     \hfill
     \begin{subfigure}[b]{\linewidth}
         \centering
        \includegraphics[width = .82
    \linewidth]{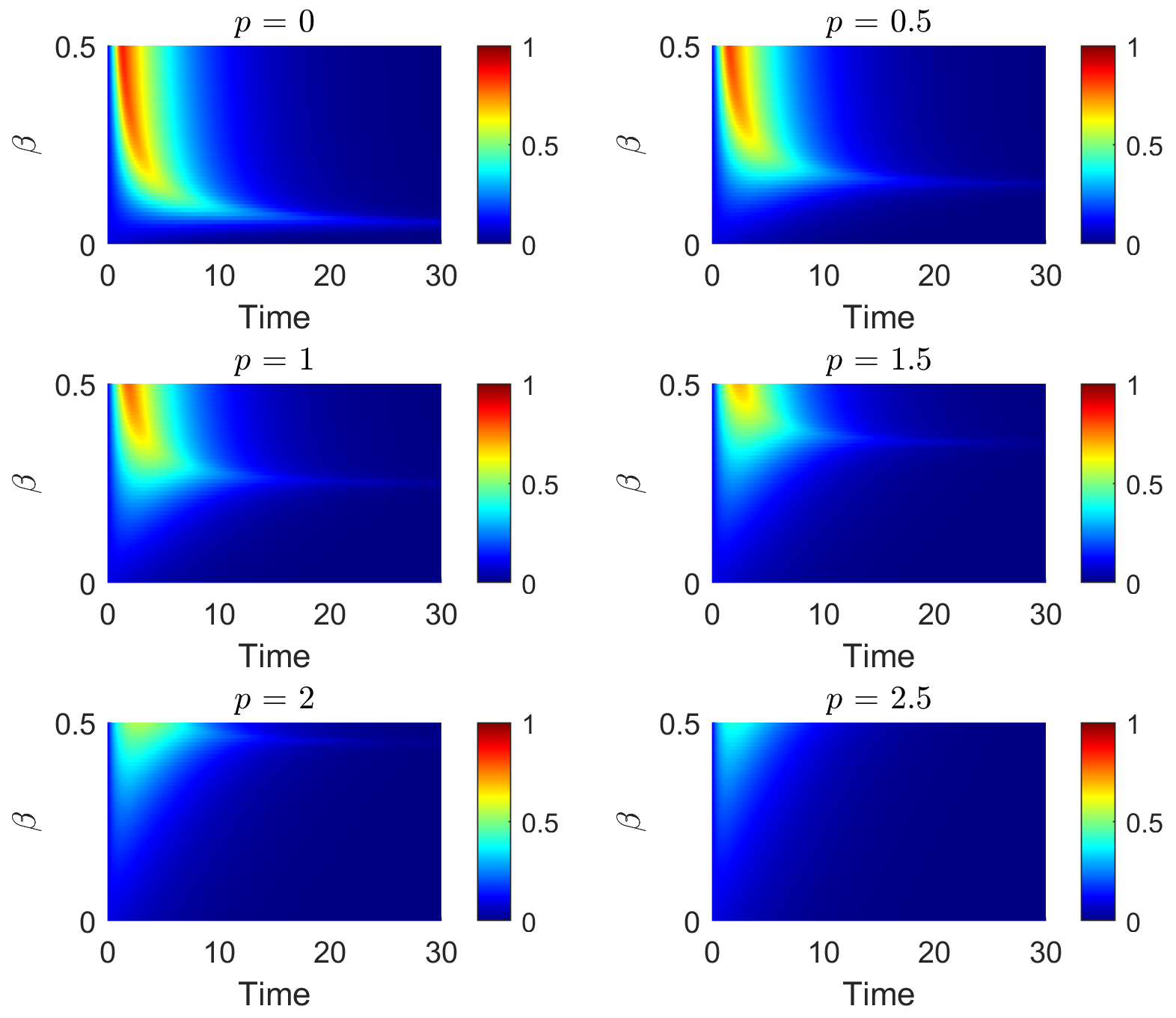}
    \caption{Compartment Model}
    \label{fig:ODEContour}
     \end{subfigure}
     \begin{subfigure}[b]{\linewidth}
         \centering
        \includegraphics[width = .82
    \linewidth]{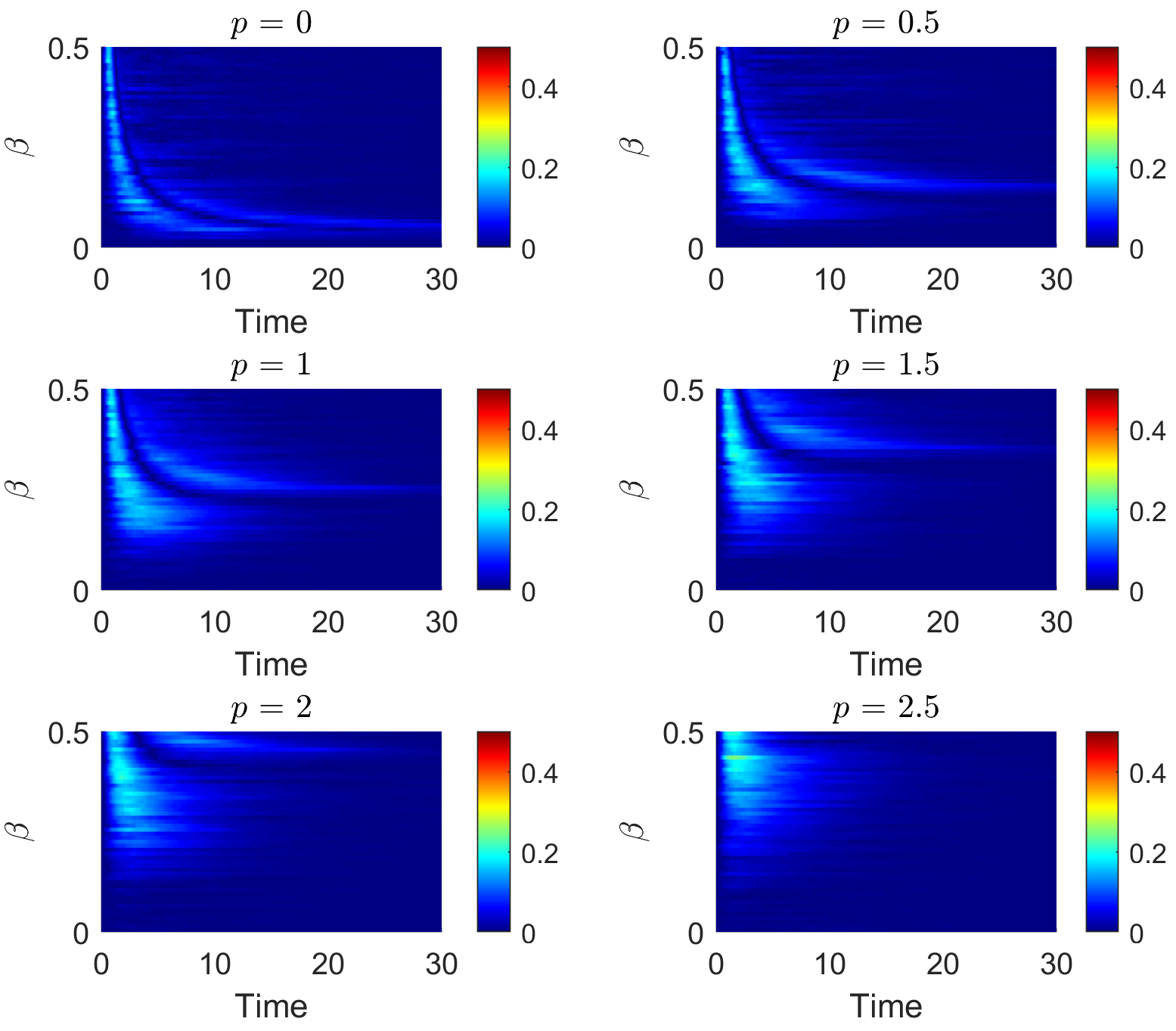}
    \caption{Absolute Difference}
    \label{fig:DifferenceContour}
     \end{subfigure}
        \caption{Plots of infected population proportions as functions of time and $\beta$ for various values of $p$ for (a) network model, (b) compartment model and (c) absolute difference between models.}
    \label{fig:Contours}
\end{figure}

The size and duration of an infectious event, as measured through the infected population, are key components to understanding the severity of an outbreak. Similarly, the infected nodes and the $I$ compartment are characteristic of overall model dynamics. In Figure \ref{fig:Contours}, we plot the proportion of the population that is infected over time and $\beta$ for $p$ values of $0$, $0.5$, $1$, $1.5$, $2$, and $2.5$ in each subplot. Figure \ref{fig:Contours}(a) has results for the network model. Striations on the plot are a result of the coarseness of simulations on a 100 node network. This plot demonstrates the influence of the deactivation rate as the contours appear to shift upward, towards higher $\beta$ values, for higher $p$ values. With $p = 0$, the infected proportion remains less than $0.2$ for all time for only $\beta< 0.05$. Conversely, with $p = 2.5$, the infected proportion is less than $0.2$ for all time for $\beta <0.2$. The difference in $\beta$ values corresponding to infected populations of the same size indicates that the deactivation rate lowers the effective infection rate as we would expect. 

Figure \ref{fig:Contours}(b) repeats the same plots described above for the compartment model. These plots appear to be roughly the same as those shown for the network model indicating convergence of our compartment model to the network model. The compartment model plots also show a ``tail" where the infected population is non-zero for an extended period of time for a particular $\beta$ value in each $p$ value subplot. Below this tail, the infection has a lower peak and shorter duration than in simulations above the tail. This tail corresponds to a $\beta$ and $p$ combination in which the recovery rate is approximately balanced by the effective infection rate causing a prolonged infectious event as the infected population proportion remains roughly constant. This behavior is hard to detect in the network model given the small, finite number of nodes used in simulations. 

Finally, to justify the convergence of the network model to the compartment model, Figure \ref{fig:Contours}(c) shows the absolute difference between the results shown in Figures \ref{fig:Contours}(a) and (b). Note, the scale for the infected proportion only ranges from $0$ to $0.5$ in this plot. This figure demonstrates remarkable consistency between the infected proportions in the network and the compartment model simulations. The maximum difference in value is less than $0.2$ across all plotted values and most inconsistencies are in the the peak value of the infected proportion and along the tail seen in the compartment model.

\subsection{Comparison of Nodes and Edges over Time}

\begin{figure}[htp]
     \centering
     \begin{subfigure}[b]{\linewidth}
         \centering
         \includegraphics[width = \linewidth]{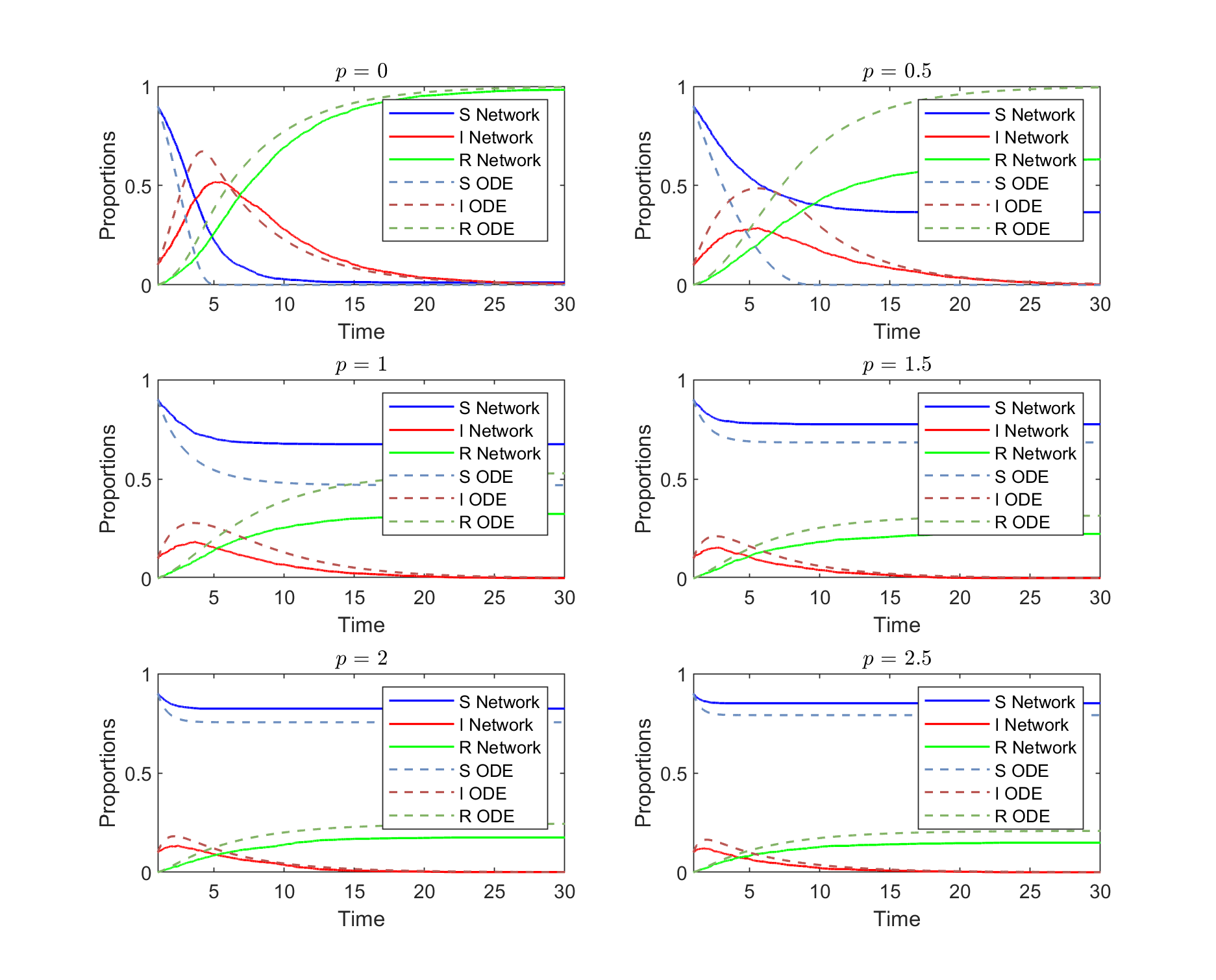}
         \caption{Node States}
         \label{fig:NodesTime}
     \end{subfigure}
     \hfill
     \begin{subfigure}[b]{\linewidth}
         \centering
        \includegraphics[width =
    \linewidth]{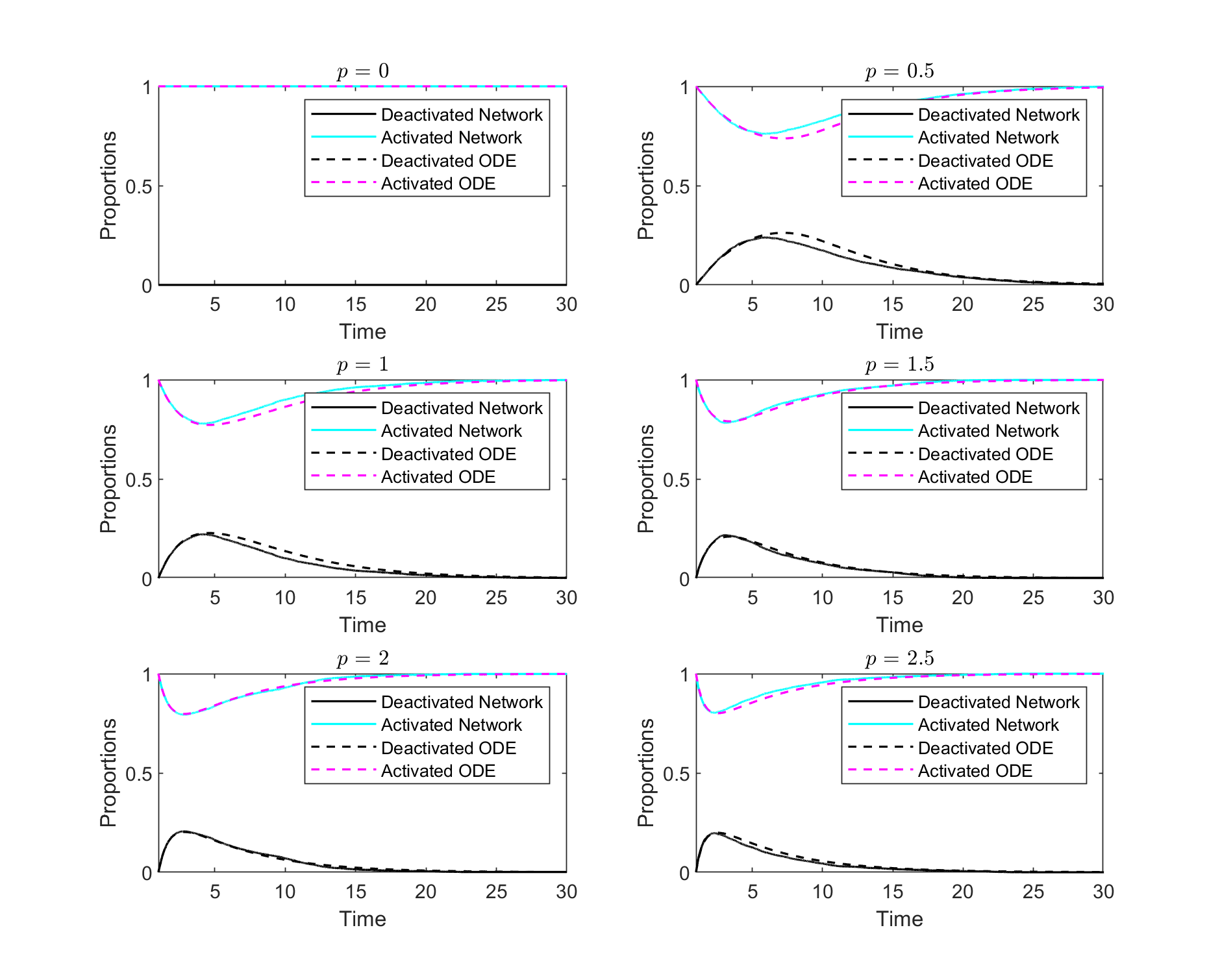}
         \caption{Edge States}
         \label{fig:EdgesTime}
     \end{subfigure}
        \caption{Average proportions in SIR model of (a) node states and (b) edge conditions over time with $\beta = 0.2$, $\gamma = 0.2$, $r = .9$, $p$ ranging from 0 to 2.5 and $I_0 = 10$. Solid lines correspond to results from the network model and dashed correspond to results from the compartment ODE model.}
    \label{fig:SIRNodesEdges}
\end{figure}

While the previous section analyzed the convergence of the infected population between our models, it is also important to consider the consistency in the other node and edge proportions between the models. Figure \ref{fig:SIRNodesEdges} plots the proportion of (a) node types and (b) edge types over time for $\beta = 0.2$ and $p \in [0, 0.5, 1, 1.5, 2, 2.5]$. Overall, Figure \ref{fig:SIRNodesEdges}(a) show consistency between the network and compartment model. The most notable discrepancy is an over estimate of the infected population by the compartment model. These plots also demonstrate the influence of the $p$ value in the proportion of the population that is in the recovered category after the disease has died off. Since an $SIR$ model assumes immunity, the ending recovered population proportion is equivalent to the cumulative proportion of the population which was infected over the course of the epidemic. This factor is significant to assessing the severity of an outbreak and will be revisited in the following subsection. 

Figure \ref{fig:SIRNodesEdges}(b) demonstrates remarkable consistency in the proportion of edges deactivated over time between the network and compartment model. This is also a key indicator that the compartment model we constructed matches the assumptions made in the network model and strengthens our use of the compartment model for broader system analysis.

\subsection{Epidemic Severity Analysis}

Traditionally, the basic reproduction rate, $R_0$, is used to assess the severity of an outbreak. Since this value quantifies the conditions under which the number of infected individuals grows in time, the rate of growth of the infected compartment can be used to calculate a formula for $R_0$. Using the $\dot{I}$ equation from Equation \eqref{eq:ODENodes}, we take the limit as the initial infected population, $I(0)=I_0$, goes to zero at time $t=0$. To compute this limit, we approximate $[SI](0)$ as $\langle k \rangle I_0$ given that it is expected all edges from the initially infected nodes will be connected to susceptible nodes as $I_0 \to 0$. This gives the calculation
\begin{equation}
     \lim_{I_0\to0} \frac{1}{I_0}\dot{I}(0) = \lim_{I_0\to0} \frac{1}{I_0} \left [\beta \langle k \rangle I_0 -\gamma I_0 \right]= \beta\langle k\rangle - \gamma. 
\end{equation}
Setting this equation equal to $0$ produces the formula
\begin{equation}
    R_0 = \frac{\beta \langle k \rangle}{\gamma}. \label{Eqn:R0}
\end{equation}

Substituting in the parameters $\langle k \rangle = 12$ and $\gamma = 0.2$ we find a critical value $\beta^* = \frac{1}{60}$ corresponding to $R_0=1$ for our simulations. For values of $\beta < \frac{1}{60}$, we expect the disease to die off immediately since the infected compartment is shrinking. For $\beta > \frac{1}{60}$, we know that the infected compartment does not decrease immediately but we cannot assume anything else about the system behavior.

From Figure \ref{fig:Contours}(b) it is clear that the severity of the disease depends on both $\beta$ and $p$ parameter values while from the above calculation $R_0$ does not. Given this limitation and the otherwise limited information provided by the standard $R_0$ calculation, we extend our analysis to the $[SI]$ edges. Since $[SI]$ is a key component of the $\dot{I}$ equation, we hypothesize that the growth of this compartment may further exemplify system behaviors. We mirror the same calculations as done above on the $\dot{I}$ equation to the $[\dot{SI}]$ equation as shown in Equation \eqref{eq:ODEEdges}. Writing out the moment closure approximation, this equation becomes

\begin{equation}
    \begin{aligned}
        \dot{[SI]} =& \beta\left(\left(\frac{\langle k \rangle - 1}{\langle k \rangle}\right) \frac{[SS][SI]-[SI]^2}{[S]}-[SI]\right)\\
        &-(\gamma+p)[SI].
\end{aligned}
\end{equation}
Note that $[SS](0)+[SI](0)=\bar{N}$ and, as used above, $[SI](0) \approx \langle k \rangle I_0$. Making the above substitutions, we have

\begin{equation}
\begin{aligned}
    \frac{\dot{[SI]}(0)}{I_0} =& \beta\left(\left(\frac{\langle k \rangle - 1}{\langle k \rangle}\right) \frac{\frac{\langle k \rangle^2N}{2}-2\langle k \rangle^2 I_0}{N-I_0}-\langle k\rangle\right)\\
    &-(\gamma+p)\langle k\rangle.
\end{aligned}
\end{equation}
Therefore,
\begin{equation} \label{eq:SIdot}
     \lim_{I_0\to0} \frac{1}{I_0}[\dot{SI}](0) = \beta \left(\frac{\langle k \rangle ^2}{2}-\frac{3}{2}\langle k \rangle\right) - (\gamma +p)\langle k \rangle.
\end{equation}
Setting this equal to $0$ and solving for $p$ gives $\dot{SI}(0)=0$ when 
\begin{equation}  \label{eq:p1beta}
     p_1^* = \beta\left(\frac{\langle k \rangle}{2}-\frac{3}{2} \right)-\gamma=\frac{\gamma}{2 \langle k \rangle }\left(R_0(\langle k \rangle -3)-2\langle k \rangle\right).
\end{equation}
Since $\ddot{S} = -\beta\dot{[SI]}$, it follows that $p_1^*$ is also the critical transition for the concavity or acceleration of the $S$ compartment. Note, necessary and sufficient conditions for $p_1^*$ to exist in the sense that it is a positive number are that $\langle k \rangle >3$ and $R_0>2\langle k \rangle /(\langle k \rangle -3)$. That is, the network has to on average have a sufficiently large number of connections and the disease has to be sufficiently contagious for deactivating edges to be necessary.

To find an equivalent critical value for the concavity of the $I$ compartment, we have 

\begin{equation}
     \lim_{I_0\to0} \frac{1}{I_0}[\ddot{I}](0) = \lim_{I_0\rightarrow 0}\left(\beta [\dot{SI}](0)-\gamma \dot{I}(0)\right).
\end{equation}
Substituting in the previously computed limit for $[\dot{SI}](0)$ as found in Equation \eqref{eq:SIdot} and solving for $p$ gives the critical value 
\begin{equation} \label{eq:p2beta}
    p_2^* = p_1^* -\gamma + \frac{\gamma^2}{\beta \langle k\rangle}=p_1^* -\gamma\left(1 - \frac{1}{R_0}\right)
\end{equation}
for the concavity or acceleration of the $I$ compartment. Note, necessary and sufficient conditions for $p_2^*$ to be a positive number are that $p_1^*>0$, which implies that $R_0>1$,  and $\gamma<p_1^*R_0/(R_0-1)$. That is, the disease has to be sufficiently contagious and the recovery rate must be not too large for deactivating edges to be necessary.

\begin{figure}
    \centering
    \includegraphics[width =
    \linewidth]{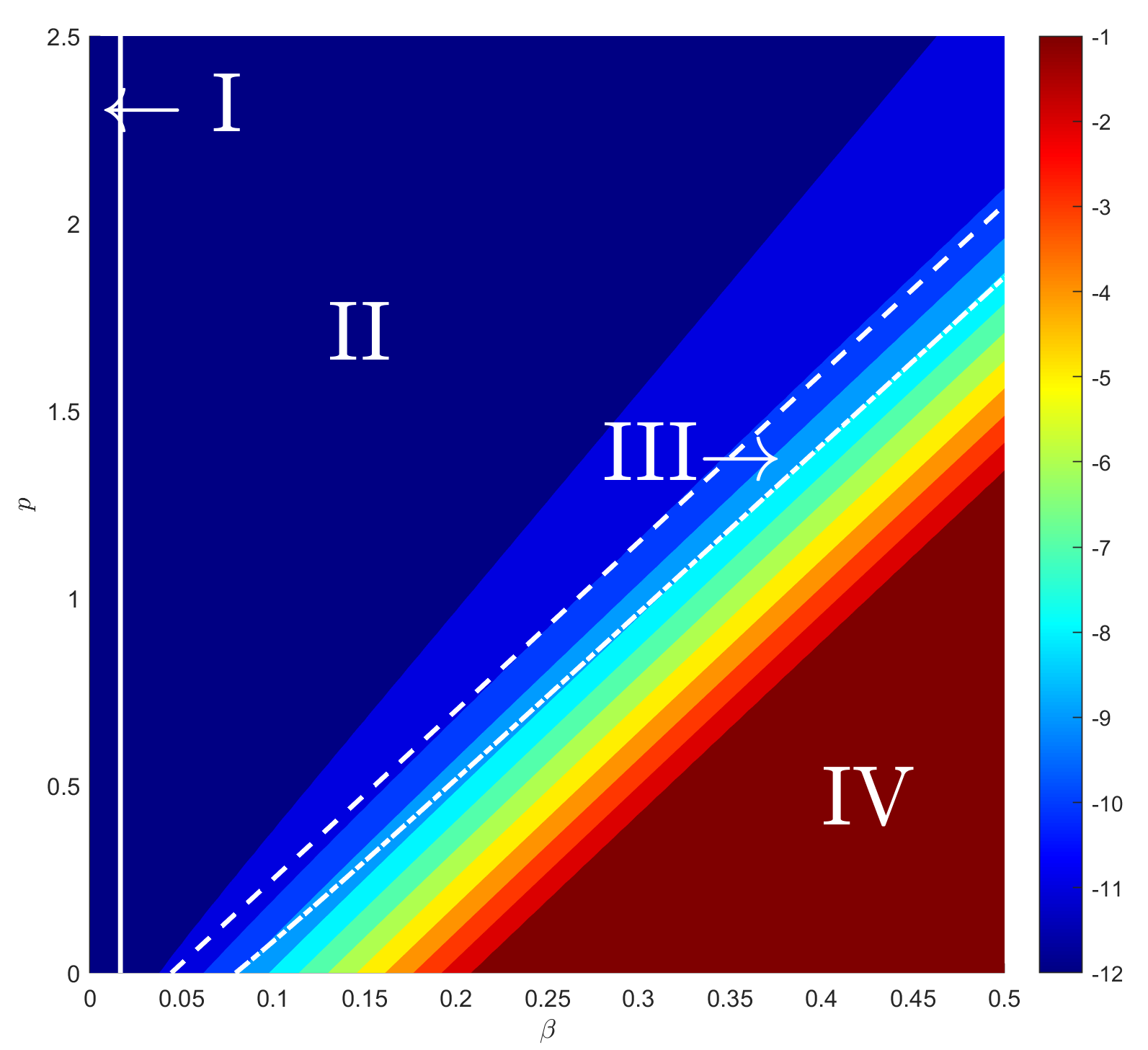}
    \caption{Contour plot of the log of the total recovered proportion of the population for the compartment model for ranging $
    \beta$ and $p$. The solid white line corresponds to $\beta=\beta^*$. The dashed white line is $p_1^*$ and the dashed-dotted line is $p_2^*$.}
    \label{fig:BetaPContour}
\end{figure}

To investigate the influence of these critical values on system behavior, we consider the cumulative proportion of the population which was infected over the course of the infectious event. In an $SIR$ model, this is equivalent to calculating the ending recovered population proportion, $R(t_f)$. Figure \ref{fig:BetaPContour} plots the log of the ending recovered population proportion for $0\leq\beta\leq 0.5$ and $0 \leq p \leq 2.5$. For these simulations, an $I_0 = 10^{-10}$ was used in correspondence with the analytic assumption of $I_0 \to 0$. The other compartments were then initialized with $S(0) = 100-I_0$, $SI(0) = \langle k \rangle I_0$, $SS(0) = \bar{N}-\langle k\rangle I_0$ and all others equal to $0$. Additionally, on Figure \ref{fig:BetaPContour},  $\beta^*$ corresponding to $R_0=1$ is plotted as a solid white line, $p_1^*$ is a dashed white line and $p_2^*$ is a dashed-dotted white line.

These critical transitions partition Figure \ref{fig:BetaPContour}  into four regions. In Region I, $R_0<1$, $\dot{I}(0)<0$, $\ddot{S}(0)<0$ and  $\ddot{I}(0)<0$. In this region the disease quickly dies out and the dynamics are equivalent to the standard $SIR$ model with $R_0<1$. In Region II, $R_0>1$, $\dot{I}(0)>0$, $\ddot{S}(0)<0$ and  $\ddot{I}(0)<0$. In this region even though $\dot{I}(0)>0$ the total number of infections is still low since the rate of change of infections is decelerating. In Region III, $R_0>1$, $\dot{I}(0)>0$, $\ddot{S}(0)>0$ and  $\ddot{I}(0)<0$. While the rate of change of infections is initially accelerating in this region, the rate of change of susceptible individuals is initially decelerating and thus again the total number of infections is still comparably low. In Region IV, $R_0>1$, $\dot{I}(0)>0$, $\ddot{S}(0)>0$ and  $\ddot{I}(0)>0$. Consequently, in Region IV the number of infections is orders of magnitude higher than in regions I-III and the dynamics is similar that of a standard $SIR$ model with $R_0>1$. It is interesting to note that the existence of Regions II and III are unique to compartment models that include edge dynamics.  


\section{Discussion}
In this work we developed and analyzed a mathematical model for the spread of an $SIR$ type infectious disease on an adaptive network with temporary link deactivation. The approach taken was to develop a system of fourteen differential equations for not only the node states but the edge states as well. This mean field approach agreed well with Monte-Carlo simulations of small sized networks. Through an analysis of these equations we not only recovered the standard calculation of $R_0$ but identified two new parameters $p_1^*$, $p_2^*$ which also control the severity of the epidemic. Specifically, while the value of $R_0$ controls whether the infection is growing in time, if the deactivation rate $p$ is below $p_2^*$ then the initial number of infections is accelerating leading to a severe epidemic, i.e. Region IV in Figure \ref{fig:BetaPContour}. This is in contrast with the standard $SIR$ model in which $R_0=1$ is the condition in which both $\dot{I}(0)$ and $\ddot{I}(0)$ changes sign. Indeed, one naive approach to understanding the dynamics caused by the deactivating of connections is to assume that deactivating edges is equivalent to lowering the value of $\beta$ in the standard $SIR$ model. Our model shows that this approach will not adequately capture the nonlinear interactions between the node and edge dynamics which are necessary to model the spread of the disease. 

It is important to note that while the dynamics of the ODE model captures the mean field dynamics of the edge states it overestimates the severity of the disease as compared to the network model. The cause for this discrepancy is at least three fold. First, the moment closure assumed that the average excess degree $\langle k \rangle_{\text{ex}}$ was equal to $\langle k \rangle -1$. However, the random variables $k$ and $k_{\text{ex}}$ have different distributions and the relationship between their averages is an inequality called the ``friendship paradox'' where $\langle k \rangle_{\text{ex}}\geq \langle k \rangle +1$ \cite{feld1991your}. In particular, in graphs in which there is a significant variance in the degree distribution, it is not clear if a set of differential equations for the various compartments can be derived in the continuum limit \cite{keeling2005networks}. Second, in the derivation of the moment closure, higher order information about the topology of the network such as clustering and the number of triangles were ignored. Third, the truncation of the system at the level of nodes and edges excludes the dynamics of higher order links which depending on the structure of the graph could be relevant. Many of these challenges can be addressed by more carefully approximating the conditional distributions that arise in the moment closure approximation; see for instance \cite{TemporaryLinkDeactivation1, TemporaryLinkDeactivation2}. Nevertheless, since the ODE models provide overestimates for the severity of the disease, the critical deactivation rates given by Equation \eqref{eqn:critica_rates} are still useful in that they provide upper bounds for the critical deactivation rates in the realized network dynamics. 

Finally, we propose that the general approach of introducing compartments for the edge dynamics discussed in this paper is the more natural approach when modeling adaptive networks. Specifically, when considering the spread of infectious diseases in which there is human behavior in the form of quarantining, contact tracing, reconnecting, etc. it is important to consider the dynamics of the connections themselves, i.e. the edges. The alternative approach of introducing new node states as compartments does not capture how the topology of the connections themselves changes during the epidemic. Indeed, this discrepancy is captured in our model due to the existence of parameters in addition to $R_0$  which depend nonlinearly on the average node degree and also govern the severity of the epidemic.


\section*{Appendix}
In this appendix we briefly derive the moment closure given by Equation \eqref{Eqn:MomentClosure} following the derivation given on page 124 of \cite{kiss2017mathematics}. First, in a network with a average node degree $\langle k \rangle$, it follows that $\langle k \rangle B$ is equal to the expected number of edges containing a node of status $B$ and thus $[AB]/(\langle k \rangle B)$ and $[BC]/(\langle k \rangle B)$ correspond to the expected proportion of edges which start at a status $B$ node that are of type $[AB]$ or $[BC]$ respectively. Therefore, if we are given that a node $B$ is connected to two other nodes then the probability that the three nodes forms a triple link of type $[ABC]$ is approximately given by $[AB][BC]/(\langle k \rangle B)^2$. Therefore, since the number of ways to choose the edges connecting to $B$ is given by $\langle k \rangle(\langle k \rangle-1)$, it follows that the probability that a triple link a node status of $B$ at its center is of type $[ABC]$ is equal to $\langle k \rangle(\langle k\rangle -1)([AB][BC])(\langle k\rangle [B])^2$. Finally, we calculate the expected value of $[ABC]$ triples by multiplying by the proportion of $B$ nodes to obtain the following moment closure approximation:
\begin{equation*}
[ABC]\approx\frac{\langle k\rangle-1}{\langle k\rangle}\frac{[AB][BC]}{B}.
\end{equation*}

\section*{Acknowledgements} 
H.S. and J.G. acknowledge the 2020 American Institute of Mathematics online summer program entitled \emph{Dynamics and data in the COVID-19 pandemic} in which some of the background research for this project was conducted. H.S. was a student in this program and J.G. an organizer. J.G. acknowledges the support of an Archie award and a Sterge Faculty Fellowship at Wake Forest University which both provided funding for this project.

\section*{Author Contributions}
H.S. and J.G. developed the theoretical formalism and performed the analytic calculations. H.S. performed the numerical simulations. All the authors discussed the results, contributed to writing and read and approved the final manuscript.







\end{document}